\documentclass[preprint]{elsarticle} % fr Rnder bei beidseitigem Druck

\usepackage{epstopdf}% To incorporate .eps illustrations using PDFLaTeX, etc.
\usepackage[caption=false]{subfig}% Support for small, `sub' figures and tables

\usepackage{setspace}

\usepackage[T1]{fontenc}        % Worttrennung mit Umlauten
\usepackage{geometry}
\usepackage[utf8x]{inputenc}    % Unicode
\usepackage{ucs}                % Unicode 
\usepackage{amsmath}            % Mathematik Erweiterungen
\usepackage{amsfonts}
\usepackage{amssymb}
\usepackage{amsthm} 
\usepackage{thmtools}
\usepackage{ifthen}
\usepackage{color}
\usepackage{graphicx} 
\usepackage{textcomp}            %Zusätzliche Symbole
\usepackage{palatino}            %Schrift
\usepackage{eulervm}            %Mathe Zeichen
\usepackage{varioref}            %Verweis von Seiten
\usepackage{numcompress}\biboptions{authoryear}
\usepackage{cleveref}             % Kluge Verweise 

\usepackage{float}
\usepackage{multirow}
\usepackage{rotating}
\usepackage{subfig}
\usepackage{xr} 
\usepackage{setspace} 
\usepackage{tabularx}
\usepackage{booktabs}
\usepackage{hhline}
\newcommand{\R}{{\mathbb{R}}}         % \R f�r R (reelle Zahlen\right)
\newcommand{\E}{{\mathbb{E}}}
 \newcommand{\N}{{\mathbb{N}}}         % \N f�r N (nat�rliche Zahlen\right)
\newcommand{\rev}[1]{\textcolor{black}{#1}}

\newcommand{\bay}{\begin{array}}
\newcommand{\eay}{\end{array}}

\newcommand{\bqa}{\begin{eqnarray*}}
\newcommand{\eqa}{\end{eqnarray*}}

\newcommand{\bee}{\begin{eqnarray*}}
\newcommand{\eee}{\end{eqnarray*}}

\newcommand{\bea}{\begin{eqnarray*}}
\newcommand{\eea}{\end{eqnarray*}}

\newcommand{\bqan}{\begin{eqnarray}}
\newcommand{\eqan}{\end{eqnarray}}

\newcommand{\be}{\begin{eqnarray}}
\newcommand{\ee}{\end{eqnarray}}

\newcommand{\bit}{\begin{itemize}}
\newcommand{\eit}{\end{itemize}}

\newcommand{\ben}{\begin{enumerate}}
\newcommand{\een}{\end{enumerate}}

\newcommand{\beq}{\begin{equation}}
\newcommand{\eeq}{\end{equation}}

\newcommand{\bdes}{\begin{description}}
\newcommand{\edes}{\end{description}}

\newcommand{\btb}{\begin{tabular}}
\newcommand{\etb}{\end{tabular}}

\newcommand{\bcen}{\begin{center}}
\newcommand{\ecen}{\end{center}}

\newcommand{\bmp}{\begin{minipage}}
\newcommand{\emp}{\end{minipage}}

\newcommand{\Cov}{\operatorname{{\it Cov}}}

\newcommand{\tr}{\operatorname{tr}}

\newcommand{\rank}{\operatorname{\it rank}}

\newcommand{\vb}{\boldsymbol{b}}

\newcommand{\vh}{\boldsymbol{h}}

\newcommand{\vp}{\boldsymbol{p}}
\newcommand{\vq}{\boldsymbol{q}}

\newcommand{\vv}{\boldsymbol{v}}

\newcommand{\vx}{\boldsymbol{x}}
\newcommand{\vy}{\boldsymbol{y}}
\newcommand{\vz}{\boldsymbol{z}}

\newcommand{\vA}{\boldsymbol{A}}

\newcommand{\vH}{\boldsymbol{H}}
\newcommand{\vI}{\boldsymbol{I}}
\newcommand{\vJ}{\boldsymbol{J}}

\newcommand{\vP}{\boldsymbol{P}}

\newcommand{\vR}{\boldsymbol{R}}

\newcommand{\vT}{\boldsymbol{T}}

\newcommand{\vV}{\boldsymbol{V}}

\newcommand{\vX}{\boldsymbol{X}}
\newcommand{\vY}{\boldsymbol{Y}}

\newcommand{\vbeta}{\boldsymbol{\beta}}

\newcommand{\vmu}{\boldsymbol{\mu}}

\newcommand{\vSigma}{\boldsymbol{\Sigma}}

\newcommand{\vtheta}{\boldsymbol{\theta}}

\newcommand{\veins}{{\bf 1}}
\newcommand{\vnull}{{\bf 0}}

\DeclareMathOperator{\vech}{vech}

\newtheoremstyle{Test1}% name of the style to be used
  {2 \baselineskip}% measure of space to leave above the theorem. E.g.: 3pt
  {1.5 \baselineskip}% measure of space to leave below the theorem. E.g.: 3pt
  {\itshape}% name of font to use in the body of the theorem
  {-0.0ex}% measure of space to indent
  {\fontfamily{ppl}\fontseries{l}\fontshape{n}}% name of head font
  {:}% punctuation between head and body
  {\newline}% space after theorem head; " " = normal interword space
   {}% Manually specify head

\theoremstyle{Test1}
\newtheorem{Sa}{Theorem}[section]
\newtheorem{theorem}{Theorem}[section]
\newtheorem{re}[Sa]{Remark}

\newtheorem{Ko}[Sa]{Corollary}

\newcolumntype{x}[1]{!{\centering\arraybackslash\vrule width #1}}

\makeatletter
\renewenvironment{proof}[1][\proofname]{\par
  \pushQED{\qed}%
  \fontfamily{ppl}\fontseries{m}\fontshape{it} \topsep6\p@\@plus6\p@\relax
  \trivlist
  \item[\hskip\labelsep
        \bfseries
    #1\@addpunct{:}]\ignorespaces
}{%
  \popQED\endtrivlist\@endpefalse
}
\makeatother

\begin{document}

\title{Choice of the hypothesis matrix for using the Wald-type-statistic}
\author[1]{Paavo Sattler\corref{cor1}%
}%\fnref{fn1}
\ead{paavo.sattler@tu-dortmund.de}
\affiliation[1]{Institute for Mathematical Statistics and Industrial Applications, Faculty of Statistics, Technical
University of Dortmund, Joseph-von-Fraunhofer-Strasse 2-4, 44221 Dortmund, Germany}
\author[2]{Georg Zimmermann}
\affiliation[2]{Team Biostatistics and Big Medical Data, IDA Lab Salzburg, Paracelsus Medical University Salzburg, Strubergasse 16, 5020 Salzburg, Austria}
\cortext[cor1]{Corresponding author}

\begin{abstract}
A widely used formulation for null hypotheses in the analysis of multivariate $d$-dimensional data is $\mathcal{H}_0: \vH \vtheta =\vy$ with $\vH\in\R^{m\times d}$, $\vtheta\in \R^d$ and $\vy\in\R^m$, where $m\leq d$. Here the unknown parameter vector $\vtheta$ can, for example, be the expectation vector $\vmu$, a vector $\vbeta $ containing regression coefficients or a quantile vector $\vq$. Also, the vector of nonparametric relative effects $\vp$ or an upper triangular vectorized covariance matrix $\vv$ are useful choices.
However, even without multiplying the hypothesis with a scalar $\gamma\neq 0$, there is a multitude of possibilities to formulate the same null hypothesis with different hypothesis matrices $\vH$ and corresponding vectors $\vy$. Although it is a well-known fact that in case of $\vy=\vnull$ there exists a unique projection matrix $\vP$ with $\vH\vtheta=\vnull\Leftrightarrow \vP\vtheta=\vnull$, for $\vy\neq \vnull$ such a projection matrix does not necessarily exist.
\noindent
Moreover, since such hypotheses are often investigated using a quadratic form as the test statistic, the corresponding projection matrices often contain zero rows; so, they are not even effective from a computational aspect. In this manuscript, we show that for the Wald-type-statistic (WTS), which is one of the most frequently used quadratic forms, the choice of the concrete hypothesis matrix does not affect the test decision. Moreover, some simulations are conducted to investigate the possible influence of the hypothesis matrix on the computation time.\end{abstract}

\maketitle

Formulating an appropriate null hypothesis is an important aspect of statistical data analysis. For multivariate analysis of an unknown $d$-dimensional parameter vector $\vtheta$, a commonly used hypothesis is given as   $\mathcal{H}_0: \vH \vtheta =\vy$. Here $\vH\in\R^{m\times d}$ denotes the so-called hypothesis matrix, and  $\vy\in\R^m$, where $m\leq d$. 
Thereby, the parameter vector $\vtheta$ depends on the model and the underlying research question of interest. For example, $\vtheta$ may consist of expectation vectors of one or more groups of observations in a one- or multi-way analysis of variance ({\it e.g.}, \cite{brunner2017} and \cite{friedrich2017permuting}).
In \cite{friedrich2016} and \cite{rubarth2022}\rev{,} the nonparametric relative effect $\vp$, which is an important effect measure in nonparametric statistics, is used. Furthermore, also vectors  of quantiles $\vq$ or vectorized covariance matrices $\vv$ might be considered ({\it e.g.,} \cite{ditzhaus2021} and \cite{sattler2022}). In order to demonstrate that these frequently encountered settings are within the scope of the above-mentioned general null hypothesis formulation, we consider the following three examples in more detail now:\\
A): Let $X_{i1},\ldots,X_{in_i}$, $i \in \{1,2,3\}$ denote three independent samples of \rev{i.i.d.} random variables with expectations $\mu_i := \E(\vX_{i1})$. In such a classical one-way analysis-of-variance setting with \rev{three} groups, it is usually of interest to compare the expectations $\mu_1, \mu_2, \mu_3$, that is, to consider the null hypothesis $\mu_1 = \mu_2 = \mu_3$. Now, stacking together the three expectations into one vector, $\vmu=(\mu_1,\mu_2,\mu_3)^\top$, this null hypothesis  could be formulated with $\vH\vmu=\vnull$, where each of the following three matrices might be plugged in for $\vH$ (actually, there are many more possible choices): 
\[\vH_1=\begin{pmatrix}
1 &-1 &0\\
0 & 1 &-1\\
1 & 0 &-1
\end{pmatrix},
 \quad \vH_2=1/3\begin{pmatrix}
2 &-1 &-1\\
-1 & 2 &-1\\
-1 & -1 &2
\end{pmatrix},
\quad\vH_3=\begin{pmatrix}
1 &-1 &0\\
0 & 1 &-1\\

\end{pmatrix}.\]

B): Consider the same situation as above, yet assuming that the random variables are ordinal, with possibly different distribution functions. Consequently, expectations cannot be used any more. Instead, for example, an important effect measure could be the so-called pairwise relative effect, which is given through $p_{ij}:=P(X_{i1}<X_{j1})+1/2P(X_{i1}=X_{j1})$ for  $i,j\in\{1,2,3\}$. One of the frequently used null hypotheses here is $p_{12}=p_{23}=p_{31}=1/2$, which means that no group tends to have smaller or bigger values than the other one.
Possible hypothesis matrices for the parameter vector $\vp=(p_{12},p_{23},p_{31})^\top$ would be 

\[\vH_1=\begin{pmatrix}
1 &0 &0\\
0 & 1 &0\\
0 & 0 &1
\end{pmatrix}  \quad \text{and}\quad \vH_2=\begin{pmatrix}
1 &-1 &0\\
0 & 1 &-1\\
0 & 0 &1
\end{pmatrix}.\]
\\
These correspond to the formulations $\vH_1\vp = (1/2,1/2,1/2)^\top$ and $\vH_2\vp = (0,0,1/2)^\top$, respectively.\\
C): Now, we consider a 2-dimensional random vector $\vX$ with covariance matrix $\vV$. If we are interested in testing whether this matrix is the 2-dimensional identity, this could be done by means of the vectorized covariance matrix $\vv:=(v_{11},v_{12},v_{22})^\top$. Thereby, $v_{rs}$ denotes the element in row $r$ and column $s$ of $\vV$. Suitable hypothesis matrices would be

\[\vH_1=\begin{pmatrix}
1 &0 &0\\
0 & 1 &0\\
0 & 0 &1
\end{pmatrix}  \quad \text{and}\quad \vH_2=\begin{pmatrix}
1 &0 &0\\
0 & 1 &0\\
1 & 0 &-1
\end{pmatrix}\]
and corresponding vectors $\vy_1=(1,0,1)^\top$ and $\vy_2=(1,0,0)^\top$, which enables us to formulate the null hypothesis (\textit{i.e.,} that the covariance matrix is equal to the 2-dimensional identity matrix) as $\vH_1\vv = \vy_1$ or $\vH_2\vv = \vy_2$, respectively.\\

These few examples show the wide range of possible parameters and corresponding hypotheses of this class of null hypotheses. However, we have also observed that each of the above-mentioned null hypotheses can be formulated with a multitude of different hypothesis matrices. These matrices might have different numbers of rows or different numbers of non-zero elements; moreover,  some of these matrices are idempotent and symmetric, while others are neither. Thus, the question arises which hypothesis matrix to choose, and importantly, whether the choice of the hypothesis matrix has an impact on the results of the subsequent statistical analysis (i.e., hypothesis tests). 

For $\vy=\vnull$ (\textit{i.e.,} $\vH\vtheta = \vnull$), there is a general consensus to use the projection matrix given through $\vP=\vH^\top(\vH\vH^\top)^+\vH$, where $\vA^+$ denotes the Moore-Penrose-Inverse of a matrix  $\vA$. An important advantage of this matrix is its uniqueness. So if $\vH_1\vtheta=\vnull$ and $\vH_2\vtheta=\vnull$ are two ways to express the same hypothesis, this hypothesis can also be formulated by $\vP\vtheta=\vH_1^\top(\vH_1\vH_1^\top)^+\vH_1\vtheta=\vH_2^\top(\vH_2\vH_2^\top)^+\vH_2\vtheta=\vnull$ (see, e.g. \cite{brunnerPuri}). Thus, the appropriate hypothesis matrix choice is straightforward.Unfortunately, for $\vy\neq \vnull$,  the solution sets of $\vH\vtheta=\vy$ and $\vP\vtheta=\vy$ are not the same. Possibly, by using an appropriate vector $\widetilde \vy$, \rev{the} equivalence of the hypotheses $\vH\vtheta=\vy$ and $\vP\vtheta=\widetilde\vy$ might be obtained, but no universally applicable technique for constructing such $\widetilde\vy$ exists. Another solution could be \rev{shifting the original data} to reformulate the hypothesis with $\vy = \vnull$. However, this might only work for some parameter $\vtheta$, such as expectations or quantiles, yet not for others like the vectorized covariance matrix or relative effects, because the latter ones are invariant under such shifts. Hence, this approach to solving the problem is not universally applicable either. Moreover, even for $\vy=\vnull$, the unique projection matrix $\vP$ is always quadratic and often seems inefficient and unwieldy. This can be seen in example A). The hypothesis matrix $\vH_2$ is the unique projection matrix with no zero entries, rational numbers as entries and three rows.
So in comparison with both other matrices, it is less attractive.

For both reasons, we will show that the value of a popular quadratic form, the so-called Wald-type-statistic (WTS), does not depend on the chosen hypothesis matrix. 
\section*{Main result}
Let $\vX$ be the available data set, which can be from one or more groups, and $N$ \rev{is} the total sample size.

Then, the above mentioned WTS is based on \rev{an} appropriate vector valued statistic $\vT(\vX)$, given through \[WTS(\vH,\vy)=N\cdot(\vH\vT(\vX)-\vy)^\top  (\vH \vSigma \vH^\top)^+(\vH\vT(\vX)-\vy)\]
with $\vSigma=\Cov(\vT(\vX))\geq0$, see for example \cite{munzel1999}. This vector $\vT(\vX)$ usually converges to a multivariate normal distribution and, in investigating expectation values, would possibly consist of one or more weighted mean vectors. In the case of asymptotic normality of $\vT(\vX)$, the limit distribution of the WTS is pivot (see, e.g. \cite{rao}), which makes it attractive for permutation approaches and similar techniques, see for example \cite{pBK} and \cite{ditzhaus2020b}.
Now we investigate the conjunction between hypothesis matrices and the corresponding projection matrices.

\begin{theorem}\label{theo1}
Let $\vH_1 \vx =\vy_1$ and $\vH_2 \vx =\vy_2$ be two systems of linear equations with the same non-trivial solution set, while $\vH_1\in\R^{m_1\times d}, \vH_2\in\R^{m_2\times d},  \vy_1\in\R^{m_1},\vy_2\in\R^{m_2}$ and $m_1,m_2\leq d$. Then, it follows that $ \vP_1=\vP_2$ holds.
\end{theorem}
This result is already known for $\vy=\vnull$ and quadratic matrices, where even an ``if and only if'' relation exists. But, this substantially more general result implies that the corresponding Wald-type statistics are the same, too.   

\begin{Ko}\label{Ko1}
Let $\vH_1 \vx =\vy_1$ and $\vH_2 \vx =\vy_2$ be two systems of linear equations with the same non-trivial solution set using $\vH_1\in\R^{m_1\times d},\vH_2\in\R^{m_2\times d},  \vy_1\in\R^{m_1}, \vy_2\in\R^{m_2}$ and $m_1,m_2\leq d$. Then  $\vH_1 \vtheta =\vy_1$ and $\vH_2 \vtheta =\vy_2$ are two ways to express the same null hypothesis, and for $\vSigma>0$ the values of the two corresponding Wald-type-statistics are the same.

\end{Ko}

This result thus ensures that the choice of the hypothesis matrix does not influence the test result as long as the chosen matrix appropriately reflects the hypothesis under consideration.

\begin{re}\phantom{1}\vspace*{-0.5cm}
\begin{itemize}

\item Under similar conditions, the result also holds for the so-called modified ANOVA-type-statistic (MATS), which was introduced in \cite{friedrich2017mats} and has a similar construction as the WTS (see more details in the supplementary material).
\item This result does not mean that for  $\vH \vtheta =\vy$ an equivalent null hypothesis can be formulated by using $\vP$, which is only the case for $\vy=\vnull$. 
\item For any hypothesis $\vH\vtheta=\vy$, alternative systems of linear equations can be developed by using elementary row operations and removing rows only containing zeros.
%\item The condition $\vSigma>0$ can be relaxed at the price of imposing additional model assumptions of the data $\vX$.

\end{itemize}
\end{re}
Since the chosen hypothesis matrix does not affect the WTS, \rev{using a simple and numerically efficient matrix is sensible}. Even for $\vy=\vnull$, the unique projection matrix $\vP$ is frequently not the best choice as it is often relatively expensive from a numerical point of view, and somewhat complicated with respect to interpretability.
From a computational view, choosing a matrix $\vH$ with $m=\rank(\vH)$ and, therefore, a minimal number of rows is favourable. Depending on the hypothesis, this can result in an essential reduction of computation time, which is investigated in the following section.

\section*{Simulation}
In this section, we will investigate the influence of the chosen hypothesis matrix on the computation time of the WTS. Here we focus on \rev{the} following two settings:
\begin{itemize}
\item[A)]Consider two independent groups with  observations $\vX_{ik}$, $i \in \{1,2\}$, $k \in \{1,...,n_i\}$ that are identically distributed within these groups, and $\E(\vX_{ij})=\vmu_i\in\R^d$, $\Cov(\vX_{ij})=\vV_i>0$. For example, consider the case of repeated measures, where the components of $\vmu_i$ represent the time-point-specific expectations in group $i$. In this context, a frequently used null hypothesis is the hypothesis of no group effect for $\vmu=(\vmu_1^\top,\vmu_2^\top)^\top$, given through $\overline \mu_{1\cdot}=\overline \mu_{2\cdot}$.

In this setting, the projection matrix can be written as $\vP_2=\vI_2-\vJ_2/2$, where $\vI_2$ denotes the identity matrix, and $\vJ_2$ is a $2\times 2$ matrix containing only 1's. 
Then, the hypothesis of no group effect can be formulated through the hypothesis matrix $\vH_1^A=(\vP_2 \otimes \vJ_d)$ by $\vH_1^A\vmu=\vnull_{2d}$. Another way to formulate the same hypothesis would be $\vH_2^A=(1,....,1,-1,....,-1)\in\R^{1\times 2d}$ and $\vH_2^A\vmu=0$.
\item[B)] Consider a positive definite symmetric matrix $\vV\in\R^{p\times p}$ as covariance matrix of observations $\vX_1,...,\vX_{n}$. The hypothesis  $\mathcal{H}_0: \tr(\vV)=\gamma$ can be investigated by using the upper triangular vectorisation $\vech(\vV)=(v_{11},...,v_{1p},v_{22},...,v_{2p},...,v_{pp}) \in\R^{d}$, $d=p(p+1)/2$ as it was introduced in \cite{sattler2022}.
Thereby we use the vector $\vh_p:=(1,\vnull_{p-1}^\top, 1, \vnull_{p-2}^\top, \dots, 1,0,1)^\top\in \R^{d}$  to define two hypothesis matrices $\vH_1^B =\vh_p\cdot \vh_p^\top$ or $\vH_2^B=  \vh_p^\top\in \R^{1\times d}$. Then the hypothesis can be formulate as $\vH_1^B\vech(\vV)=\gamma \cdot \veins_d$ and $\vH_2^B\vech(\vV)=\gamma$.
\end{itemize}

Since we want to focus on the calculation of the WTS and not on the vector T(X), nor on the estimation of the covariance matrix of T(X) we will not test the above hypotheses but use it as motivation. To get a pure result for the quadratic form, we generate one random vector $\vX_{A)}\sim\mathcal{N}_{2d}(\vnull_{2d},\vV_{2d})$ resp. $\vX_{B)}\sim\mathcal{N}_{d}(\veins_{d},\vV_{d})$. Here $\vV_{d}$ denotes a d-dimensional compound symmetry matrix given through $\vV_{d}=\vI_d+\veins_d\veins_d^\top$. With $\vT(\vX)=\vX$ for each setting, the computation time of the WTS using both respective hypothesis matrices is calculated for different dimensions $d_{A)}=(5,10,20,50,100,200)$ resp.  $d_{B)}=(15,55,120,210,325,465)$. These dimensions $d_{B)}$ are based on the dimension of the observation vector $p=(5,10,15,20,25,30)$ and the clear connection between p and d through $d=p(p+1)/2$.

Since the WTS is attractive for bootstrap and especially permutation approaches, we measure the time of 5.000 calculations of the quadratic form, which is a usual number of bootstrap steps resp. permutations.

  The computations were run by means of the \textsc{R}-computing environment version 3.6.1 \cite{R}  on an Intel Xeon E5430 quad-core CPU
 running at 2.66 GHz using 16 GB DDR2 memory on a Debian GNU Linux 7.8. The required time in seconds is displayed in \Cref{tab:Zeit1}.

 \begin{table}[ht]
\centering

      \begin{tabular}{x{1.5pt}lx{1.5pt}r|r|r|r|r|rx{1.5pt}}
       \specialrule{1.5pt}{0pt}{0pt}
   \multicolumn{1}{x{1.5pt}lx{1.5pt}}{}
&\multicolumn{6}{ cx{1.5pt}}{calculation time in seconds}
\\ \specialrule{1.5pt}{0pt}{0pt}
   d   &5&10&20&50&100  &200 
      \\  \specialrule{1.5pt}{0pt}{0pt}
    $WTS(\vH_1^A,\vnull_{2d})$& 0.928 & 1.452 & 117.934 & 1976.501 & 1616.691 & 6724.439 \\ \hline
      $WTS(\vH_2^A,0)$ &  0.737 & 0.853 & 0.862 & 22.373 & 29.275 & 44.880 \\ 
       \specialrule{1.5pt}{0pt}{0pt}
       
       \multicolumn{7}{x{1.5pt} cx{1.5pt}}{}\\
\specialrule{1.5pt}{0pt}{0pt}
  % p   &5&10&15&20&25  &30 \\
   d& 15&55&120&210&325&465\\
    \specialrule{1.5pt}{0pt}{0pt}
       $WTS(\vH_1^B,\gamma\cdot\veins_d)$ &0.818 & 1.031 & 83.394 & 323.360 & 1370.995 & 9339.278  \\ \hline
       $WTS(\vH_2^B,\gamma)$ &0.718 & 0.742 & 0.777 & 23.281 & 26.594 & 59.630  \\        
      \specialrule{1.5pt}{0pt}{0pt}

      \end{tabular}

\caption{ Average computation time in seconds of \rev{test statistics} based on different hypothesis matrices with different dimensions.}
  \label{tab:Zeit1}
\end{table}

The average computation time for both settings and hypothesis matrices can be seen in {\Cref{tab:Zeit1}. For larger dimension $d\geq 20$ resp. $d\geq 120$(which means $p\geq 15$), using the smaller hypothesis matrix saves over $99\%$ of the calculation time, which makes the difference between minutes and multiple hours. Although the relative saving of computation time for smaller dimensions is considerable,  in this case, \rev{it's} only fractions of seconds.\\
The results for settings A and B are comparable because the relations between the number of rows of the quadratic matrix and of the smaller matrix are similar. Of course, the time-saving effect increases with the increasing number of \rev{removed} rows. 
\\\\
However, this simulation demonstrates only one aspect of how a smaller hypothesis matrix can reduce the required calculation time. It is widely known that calculating the empirical covariance is a substantial part of the needed computations, especially for higher dimensions. Depending on the setting, it could simplify the calculation to multiply the underlying data with the hypothesis matrix first and afterwards estimate the empirical covariance matrix of these transformed observations.
Similar approaches \rev{exist} for bootstrap approaches, where for example, it is often more efficient to generate $\vY\sim\mathcal{N}_m(\vnull_m,\vH^\top \vV\vH)$, compared to $\vX\sim\mathcal{N}_d(\vnull_d, \vV)$ and multiplying them with the hypothesis matrix afterwards.
Therefore, choosing a smaller hypothesis matrix for the WTS has a huge potential, in particular for large dimensions and resampling approaches.

\section*{Other quadratic forms}
Unfortunately, our result only holds for the WTS and the MATS, but not for the ATS, which is a widely used quadratic form-based test statistic (see, e.g. \cite{brunner:2001}). The ATS is given through 
\[ATS(\vH,\vy)=N\cdot(\vH\vT(\vX)-\vy)^\top  (\vH\vT(\vX)-\vy)\]
resp.
\[ATS(\vH,\vy)_s=N\cdot(\vH\vT(\vX)-\vy)^\top  (\vH\vT(\vX)-\vy)/\tr(\vH \vSigma \vH^\top),\]
where the second version could be seen as standardized ATS.
 Here the results strongly depend on the chosen hypothesis matrix, as can already be seen from simple examples. Nevertheless, no general statement exists for $\vy\neq \vnull$ on choosing the hypothesis matrix. 
In addition to uniqueness, which would be the most important criterion for such a matrix, also other properties would be meaningful. This means especially that this matrix should have a minimal number of rows and many zero entries to be computationally efficient. A matrix fulfilling this would be the reduced-row-echelon-form of a matrix, which is known from the Gauß-elimination algorithm. Thus, a hypothesis could be defined as the solution of a linear equation system, and afterwards, the reduced-row-echelon-form would be calculated and used for calculating the test statistic. Although there are different ways to define the linear equation system, for each solution set and the corresponding null hypothesis, this procedure would result in one and the same hypothesis matrix, \rev{eventually leading} to a unique test statistic.

Other meaningful properties of a hypothesis matrix, like symmetry or idempotence, would be desirable but usually are in conflict with computational advantages.\\\\
On the other hand, for $\vy=\vnull$, the ``default choice'' $\vP$ is unique but possibly contains multiple zero rows or rows\rev{,} which are linear combinations of other rows. If we reconsider example C) and the \rev{covariance matrix's structure, simple examples exist for both}. Considering a two-dimensional random vector, the unique projection matrices for testing diagonality resp. sphericity (i.e., diagonality with equal diagonal elements) would be

\[
\vH=\begin{pmatrix}
0&0&0\\
0&1&0\\
0&0&0

\end{pmatrix}\quad \text{resp.} \quad  \vH=\frac 1 2\cdot\begin{pmatrix}
1&0&-1\\
0&2&0\\
-1&0&1

\end{pmatrix}.\]\\
Both, the zero rows and the linear dependencies within the hypothesis matrices make them dispensable and time-consuming. Fortunately, the removal of zero rows does not affect the value of the ATS. Also, for linear dependencies, it is possible to remove multiple rows without affecting the value of the test statistic. Without loss of generality, let $\vH_{\bullet \ell_1},...,\vH_{\bullet \ell_k}$ be pairwise linearly dependent rows with $\vH_{\bullet \ell_2}=\beta_{2}\cdot\vH_{\bullet \ell_1},...,\vH_{\bullet \ell_k}=\beta_{k}\cdot\vH_{\bullet \ell_1} $ and $\beta_2,...,\beta_k \in\R\setminus 0$. Then we could remove $\vH_{\bullet \ell_2},...,\vH_{\bullet \ell_k}$, if we replace  $\vH_{\bullet \ell_1}$ by $\sqrt{1+\beta_1^2+...+\beta_k^2}\vH_{\bullet \ell_1}$, without influencing the value of the test statistic. So, the first row of them remains but is multiplied by the square root of its summed squared frequencies. It is clear that this also works if various rows have multiples, and it could be adapted for $\vy\neq \vnull$. This yield to the next corollary:

\begin{Ko}
Let $\vH\vtheta=\vy$ be \rev{an} arbitrary hypothesis with hypothesis matrix $\vH\in \R^{m\times d}$ and denoting with $\mathcal M_0(\vH)$ the set containing indices of zero rows of $\vH$. Then there are  $\ell\leq m$ pairwise linear independent rows of the matrix $\vH$, which allows to divide $\{1,...,m\}\setminus \mathcal M_0(\vH)$ in $\ell$ disjoint index sets $\mathcal{M}_1,...\mathcal{M}_{\ell}$. For these index sets, $\min(\mathcal M_i)$ denotes the smallest element of a set and $|\mathcal M_i|$ the cardinality of a set. For $k=1,...,\ell$ these sets containing elements $m_{k,1},....,m_{k,|\mathcal{M}_k|}$. So for $k=1,...,\ell$, $j=1,...,|\mathcal{M}_k|$ it exists $\beta_{k,j}\in \R\setminus 0$ with $\vH_{\bullet m_{k,j}}=\beta_{k,j} \cdot \vH_{\bullet \min(\mathcal M_k)}$.
Then we can express the same hypothesis through $\widetilde \vH\vtheta =\widetilde\vy$ \rev{using} 

\[\widetilde\vH=\begin{pmatrix}
w_1\cdot\vH_{\bullet \min(\mathcal M_1)}\\
w_2\cdot\vH_{\bullet \min(\mathcal M_2)}\\
\vdots\\
w_\ell\cdot\vH_{\bullet \min(\mathcal M_\ell)}
\end{pmatrix}\quad
\text{and} \quad \widetilde \vy=\begin{pmatrix}
\sqrt{1+\beta_{1,2}^2+...+\beta_{1,|\mathcal{M}_1|}^2}\cdot y_{ \min(\mathcal M_1)}\\
\sqrt{1+\beta_{2,2}^2+...+\beta_{2,|\mathcal{M}_2|}^2}\cdot y_{ \min(\mathcal M_2)}\\
\vdots\\
\sqrt{1+\beta_{\ell,2}^2+...+\beta_{\ell,|\mathcal{M}_\ell|}^2}\cdot y_{\min(\mathcal M_\ell)}
\end{pmatrix},\]
\rev{with $w_k=\sqrt{1+\beta_{k,2}^2+...+\beta_{k,|\mathcal{M}_k|}^2}$.}
Moreover, these hypothesis matrices \rev{fulfilling} $ATS(\vH,\vy)=ATS(\widetilde \vH,\widetilde \vy)$ and $ATS_s(\vH, \vy)=ATS_s(\widetilde \vH,\widetilde\vy).$

\end{Ko}

Both adaptations make the remaining matrix more efficient but neither idempotent nor symmetric. However, \rev{ if $\min(\mathcal M_1)<\min(\mathcal M_2)...<\min(\mathcal M_\ell)$}, starting with an unique projection matrix  \rev{lead to an adapted version which is} also  unique.

\section*{Conclusion}
The choice of the hypothesis matrix for the corresponding research question and the null hypothesis  $\mathcal{H}_0: \vH \vtheta =\vy$ for $\vy\neq \vnull$ are two topics that until now have received little attention in the literature. At the same time, the chosen matrix can potentially have a substantial impact on the test result and the required computation time. For the Wald-type statistic, our result allows for choosing an easily interpretable and computationally efficient hypothesis matrix without any risk of affecting the test result.\\
This is impossible for other quadratic form-based test statistics such as the ATS. Therefore we introduced an approach which leads to an unique and computationally advantageous hypothesis matrix. Especially for $\vy\neq \vnull$, \rev{where} no statements of choosing the hypothesis matrix \rev{exist}, this seems \rev{an} appropriate solution.\\
Even for $\vy=\vnull$ existing difficulties are discussed and simple ways are \rev{introduced} to make the hypothesis matrix computationally more efficient.
 All in all, the choice of the hypothesis matrix is \rev{an} important but neglected topic. So for the future, a general guideline would be desirable, which continues our approaches and deal with even more quadratic forms.

\section*{Acknowledgements}
The work of Paavo Sattler was funded by
Deutsche Forschungsgemeinschaft Grant/Award Number DFG PA 2409/3-2.
Georg Zimmermann gratefully acknowledges the support of the WISS 2025 project 'IDA-Lab Salzburg' (20204-WISS/225/197-2019 and 20102-F1901166-KZP)

\section*{Appendix}
\begin{proof}[Proof of \Cref{theo1}]
Both systems of linear equation have the same solution set, so it follows with  Theorem 2.2(ii) (\cite{rao}) that the following both affine subspaces are identical

\[\begin{array}{llr}
\vV_1:& \vH_1^-\vy_1+ (\vI_d-\vH_1^-\vH_1)\vz_1& \text{arbitrary}\ \vz_1\in\vR^d,\\
V_2:& \vH_2^-\vy_2+ (\vI_d-\vH_2^-\vH_2)\vz_2& \text{arbitrary}\ \vz_2\in\vR^d,\\
\end{array}\]
since $V_1$ and $V_2$ are the corresponding solutions. Here $\vH^-$ denotes a generalized inverse of the matrix $\vH$. From Theorem 2.4(b) of \cite{rao} it is known that $\vH_i^\top(\vH_i\vH_i^\top)^+$ is a g-inverse of $\vH_i$, which leads to

\[\begin{array}{llr}
V_1:& \vH_1^\top(\vH_1\vH_1^\top)^+\vy_1+ (\vI_d-\vP_1)\vz_1& \text{arbitrary}\ \vz_1\in\vR^d,\\
V_2:& \vH_2^\top(\vH_2\vH_2^\top)^+\vy_2+ (\vI_d-\vP_2)\vz_2& \text{arbitrary}\ \vz_2\in\vR^d.\\
\end{array}\]\\\\
The equality of these two affine subspaces $V_1$ and $V_2$ means that the corresponding subspaces are equal (see, e.g. \cite{snapper2014} Proposition 14.1).
Therefore we have the equality of %the 
%Es ist bekannt das bei Gleichheit von 2 affinen Unterräumen die jeweiligen Untervektorräume gleich sein müssen(\cite{snapper2014} Proposition 14.1).  Daraus ergibt sich mit $\widetilde \vP_i=\vI_d-\vP_i$ die Gleichheit der folgenden beiden Untervektorräume: 

\[\begin{array}{llr}
U_1:& (\vI_d-\vP_1)\vz_1& \text{arbitrary}\ \vz_1\in\vR^d,\\
U_2:& (\vI_d-\vP_2)\vz_2& \text{arbitrary}\ \vz_2\in\vR^d.\\
\end{array}\]
These matrices $\vI_d-\vP_1$ and $\vI_d-\vP_2$ are also projection matrices. Since projection matrices on subspaces are unique,
$\vI_d-\vP_1=\vI_d-\vP_2$ and finally $\vP_1=\vP_2$ follows.\\
\end{proof}

So, again the projection matrix $\vP=\vH^\top(\vH\vH^\top)^+\vH$ is unique, but it \rev{does} not necessarily test the initial hypothesis. However, this result allows to prove that the value and, therefore, the test decision of the WTS \rev{does} not depend on the used hypothesis matrix.

\begin{proof}[Proof of \Cref{Ko1}]
Since the solution set is non-trivial, \rev{at least one $\vb\in \R^d$ exists} with $\vH_1\vb=\vy_1$ and $\vH_2\vb=\vy_2$.
Define $\widetilde \vH_i =\vH_i\vSigma^{1/2}$  for $i=1,2$ %and $\widetilde \vT(\vX)=\vSigma^{-1/2}\vT(\vX)$
 and calculate
\[\begin{array}{ll}WTS(\vH_1,\vy_1)&=N\cdot (\vH_1\vT(\vX)-\vy_1)^\top  (\vH_1 \vSigma \vH_1^\top)^+(\vH_1\vT(\vX)-\vy_1)\\
&=N\cdot (\vT(\vX)-\vb)^\top(\vSigma^{-1/2})^\top (\vSigma^{1/2})^\top\vH_1^\top (\vH_1 \vSigma^{1/2}(\vSigma^{1/2})^\top\vH_1^\top)^+ \vH_1 \vSigma^{1/2} \vSigma^{-1/2}(\vT(\vX)-\vb)\\
&=N\cdot (\vT(\vX)-\vb)^\top(\vSigma^{-1/2})^\top \widetilde \vH_1^\top (\widetilde\vH_1 \widetilde\vH_1^\top)^+\widetilde \vH_1 \vSigma^{-1/2}(\vT(\vX)-\vb),
\end{array}\]
and analogously
\[\begin{array}{ll}WTS(\vH_2,\vy_2)&=N\cdot (\vT(\vX)-\vb)^\top(\vSigma^{-1/2})^\top \widetilde \vH_2^\top (\widetilde\vH_2 \widetilde\vH_2^\top)^+\widetilde \vH_2 \vSigma^{-1/2}(\vT(\vX)-\vb).
\end{array}\]
If now $\widetilde \vH_1 \vx =\vy_1$ and $\widetilde \vH_2 \vx =\vy_2$ are two systems of linear equations with the same non-trivial solution set we can use \Cref{theo1} and it follows $\widetilde \vH_1^\top (\widetilde\vH_1 \widetilde\vH_1^\top)^+\widetilde \vH_1=\widetilde \vH_2^\top (\widetilde\vH_2 \widetilde\vH_2^\top)^+\widetilde \vH_2$ and therefore $WTS(\vH_1,\vy_1)=WTS(\vH_2,\vy_2)$.
\\\\
First $\vSigma^{-1/2} \vb$ solves both equations since $\widetilde \vH_1 (\vSigma^{-1/2}\vb)=\vH_1\vb=\vy_1$ and $\widetilde \vH_2 (\vSigma^{-1/2}\vb)=\vH_2\vb=\vy_2$, so the solution set is not empty. Assume now there exists a $\vx'\in\R^d$ with  $\widetilde \vH_1 \vx' =\vy_1$ but $ \widetilde \vH_2 \vx' \neq\vy_2$. Then for $\vSigma^{1/2}\vx'$ it would hold
$\vH_1(\vSigma^{1/2}\vx')=\widetilde\vH_1\vx'=\vy_1$ and $\vH_2(\vSigma^{1/2}\vx')=\widetilde\vH_2\vx'\neq\vy_2$. Since the solution set of $\vH_1 \vx =\vy_1$ and $\vH_2 \vx =\vy_2$ are the same, such a $\vx'$ can not exists. Therefore $\widetilde \vH_1 \vx =\vy_1$ and $\widetilde \vH_2 \vx =\vy_2$ are two systems of linear equations with the same non-trivial solution set and the result follows.\\
\end{proof}
%In the above proof, $\vSigma>0$ and therefore the existence of $\vSigma^{-1}$ was only necessary to factor out $\vSigma^{1/2}$. Therefore, if this is ensured in another way, this condition would be 
%unnecessary. Consider for example $d$-dimensional vectors $\vX_1,....\vX_{N}$ given through $\vX_{k}=\vmu+\vSigma^{1/2}\vep_k$ with $\vep_k\stackrel{i.i.d.}\sim\mathcal N(\vnull,\vSigma)$ with $k=1,...,N$. Then with $\vT(\vX)=\frac{1}{N}\sum_{k=1}^N \vX_k$ and $\vy=\vmu$ no additional conditions on the covariance matrix $\vSigma$ would be necessary.\\\\

%This would be for example be the case for $\vT(\vX)= \overline \vX$ and $\vy=\vnull$.\\\\

The MATS, which was introduced in \cite{friedrich2017mats}, is given through
\[MATS(\vH,\vy)=(\vH\vT(\vX)-\vy)^\top  (\vH \vSigma_0 \vH^\top)^+(\vH\vT(\vX)-\vy),\]
while $\vA_0$ denotes the matrix only containing the diagonal elements of $\vA$. Here additionally, $\vSigma_0>0$ is required, which is sufficient to get \rev{the} same result also for the MATS, by only replacing  $\vSigma$ by $\vSigma_0$ in the above proof.\\
For the ATS, even the multiplication with \rev{a} non-singular diagonal matrix \rev{usually has} a substantial effect on the value of the test statistic, and, therefore, the test decision. In the following corollary, we show some adoptions which can be conducted without \rev{affecting} the value of the test statistic.

\begin{proof}
First, we consider the equality of both hypotheses and, therefore, the linear equation systems. So assume there exists a $\vtheta' \in \R^d$ with $\vH\vtheta'=\vy$ but $\widetilde \vH \vtheta'\neq \widetilde\vy$. Then there exists at least one $k\in \{1,...,\ell\}$ with $\widetilde \vH_{\bullet k}\vtheta'\neq \widetilde y_k$. But from the definition of $\widetilde \vH$ and $\widetilde\vy$ we get
\[\begin{array}{ll}
\widetilde \vH_{\bullet k}\vtheta'
&=w_k\cdot\vH_{\bullet \min(\mathcal M_k)}\vtheta'\\
&=w_k\cdot\vy_{\min(\mathcal M_k)}\\
&= \widetilde y_k.
\end{array}\]
This is a contradiction, and therefore there can exist no $\vtheta'$ like this.\\\\
For the other direction, we first mention that for $k=1,...,\ell$, $j=1,...,|\mathcal{M}_k|$ from $\vH_{\bullet m_{k,j}}=\beta_{k,j} \cdot \vH_{\bullet \min(\mathcal M_k)}$ and $\vH\vtheta =\vy$ it directly follows 
$y_{ m_{k,j}}=\beta_{k,j} \cdot y_{ \min(\mathcal M_k)}$.\\
Now let us assume that there exists a $\vtheta' \in \R^d$  with $\widetilde \vH\vtheta'=\widetilde \vy$ but $ \vH \vtheta'\neq \vy$. Then there exists at least one $i\in \{1,...,m\}$ with $ \vH_{\bullet i}\vtheta'\neq y_i$. From the definition of $\mathcal M_0(\vH)$ it is clear that  $i\in\mathcal M_0(\vH)$ implies $y_i=0$ and that  $\vH_{\bullet i}$ is a zero row. Therefore $i\notin\mathcal M_0(\vH)$ and it exists $k\in \{1,...,\ell\}$ and $j \in \{1,...,|\mathcal{M}_k|\}$ with $i=m_{k,j} \in \mathcal M_k(\vH)$.

 \rev{Then from the definition of  $\widetilde \vH_{\bullet k }$ and $\widetilde y_k$ it is clear}

\[\widetilde \vH_{\bullet k }\vtheta'=\widetilde y_k\quad  \Leftrightarrow \quad w_k\cdot\vH_{\bullet \min(\mathcal M_k)}\vtheta'= w_k\cdot y_{ \min(\mathcal M_k)}.\]
Since $w_k>0$ and $\beta_{k,j}>0$, multiplication and division of these factors sustain the equation, we get
\[\begin{array}{lrl}

&w_k\cdot\vH_{\bullet \min(\mathcal M_k)}\vtheta'&= w_k\cdot y_{ \min(\mathcal M_k)}\\
\Leftrightarrow&\vH_{\bullet \min(\mathcal M_k)}\vtheta'&= y_{ \min(\mathcal M_k)}\\
\Leftrightarrow& \beta_{k,j}\cdot\vH_{\bullet \min(\mathcal M_k)}\vtheta'&=\beta_{k,j}\cdot y_{ \min(\mathcal M_k)}\\
\end{array}\]

With the definition of $\beta_{k,j}$ and $m_{k,j}$ we therefore receive
\[\begin{array}{lrl}
&\beta_{k,j}\cdot\vH_{\bullet \min(\mathcal M_k)}\vtheta'&=\beta_{k,j}\cdot y_{ \min(\mathcal M_k)}\\
\Leftrightarrow& \vH_{\bullet m_{k,j}}\vtheta'&= y_{ m_{k,j}}\\
\Leftrightarrow& \vH_{\bullet i }\vtheta'&= y_{ i}\\
\end{array}\]
which is again a contradiction. So we know that both linear equation systems are equal. For the equality of the \rev{test statistics}, we first define $\N_d=\{1,...,d\}$. Then \rev{with the rules for matrix calculation, we get} 
\[\begin{array}{ll}
&ATS(\vH,\vy)\\=&N\cdot (\vH\vT(\vX)-\vy)^\top(\vH\vT(\vX)-\vy)\\
=& N\cdot \sum_{i\in \N_m} [(\vH\vT(\vX))_{ i}-y_i]^2\\
=& N\cdot \sum_{i\in \N_m} [\vH_{\bullet i} \vT(\vX)-y_i]^2
\end{array}.\]
\rev{Now, we can decompose $\N_m$ in $\N_m\setminus \mathcal M_0$ and $M_0$, where for $i\in M_0$ again $\vH_{\bullet i}$ is a zero row and $y_i=0$. This leads to}
\[\begin{array}{ll}
 &N\cdot \sum_{i\in \N_m} [\vH_{\bullet i} \vT(\vX)-y_i]^2\\
=& N\cdot \sum_{i\in \N_m\setminus \mathcal M_0(\vH)} [\vH_{\bullet i} \vT(\vX)-y_i]^2+ N\cdot \sum_{i\in  \mathcal M_0(\vH)} [\vH_{\bullet i} \vT(\vX)-y_i]^2\\
=& N\cdot \sum_{i\in \N_m\setminus \mathcal M_0(\vH)} [\vH_{\bullet i} \vT(\vX)-y_i]^2+ 0.\\
\end{array}\]
\rev{A further disjunction of $\N_m\setminus \mathcal M_0$ in $M_1,...,M_\ell$ and the definition of $\beta_{k,j}$ leads to}
\[\begin{array}{ll}
& N\cdot \sum_{i\in \N_m\setminus \mathcal M_0(\vH)} [\vH_{\bullet i} \vT(\vX)-y_i]^2\\
=&N\cdot\sum_{k\in \N_\ell}\left([\vH_{\bullet m_{k,1}} \vT(\vX)-y_{m_{k,1}}]^2+...+[\vH_{\bullet m_{k,|\mathcal{M}_k|}} \vT(\vX)-y_{m_{k,|\mathcal{M}_k|}}]^2\right)\\
=&N\cdot\sum_{k\in \N_\ell}\left([\vH_{\bullet \min(\mathcal M_k)} \vT(\vX)-y_{\min(\mathcal M_k)}]^2+...+\beta_{k,|\mathcal{M}_k|}^2[\vH_{\bullet \min(\mathcal M_k)} \vT(\vX)-y_{\min(\mathcal M_k)}]^2\right)\\
=&N\cdot\sum_{k\in \N_\ell} \left[\sqrt{(1+\beta_{k,2}^2+...\beta_{k,|\mathcal{M}_{k}|}^2)}\cdot(\vH_{\bullet \min(\mathcal M_k)} \vT(\vX)-y_{\min(\mathcal M_k)})\right]^2.\\
\end{array}\]
\rev{Finally with the definition of $w_k$, $\widetilde\vH$ and $\widetilde \vy$ we get}
\[\begin{array}{ll}
&N\cdot\sum_{k\in \N_\ell} \left[\sqrt{(1+\beta_{k,2}^2+...\beta_{k,|\mathcal{M}_{k}|}^2)}\cdot(\vH_{\bullet \min(\mathcal M_k)} \vT(\vX)-y_{\min(\mathcal M_k)})\right]^2\\
=&N\cdot\sum_{k\in \N_\ell} \left[w_k\cdot(\vH_{\bullet \min(\mathcal M_k)} \vT(\vX)-y_{\min(\mathcal M_k)})\right]^2\\
=&N\cdot (\widetilde \vH\vT(\vX)-\widetilde \vy)^\top(\widetilde \vH\vT(\vX)-\widetilde \vy)\\
=&ATS(\widetilde \vH,\widetilde \vy)

\end{array}\]
\rev{and therefore $ATS(\vH,\vy)=ATS(\widetilde \vH,\widetilde \vy)$.}
\\

For the $ATS_s$, we consider the corresponding trace and use the same steps. Since $\vSigma$ is positive semi-definite, this leads to\\\\
$\begin{array}{ll}
\tr(\vH \vSigma\vH^\top)&=\sum_{i \in \N_m} (\vH \vSigma\vH^\top)_{ii}\\
&=\sum_{i \in \N_m} (\vH \vSigma^{1/2}(\vH\vSigma^{1/2})^\top)_{ii}\\
&=\sum_{i \in \N_m\setminus \mathcal M_0(\vH)} ||(\vH\vSigma^{1/2})_{\bullet i}||^2+\sum_{i\in \mathcal M_0(\vH)} ||(\vH\vSigma^{1/2})_{\bullet i}||^2\\
&=\sum_{i \in \N_m\setminus \mathcal M_0(\vH)} ||(\vH\vSigma^{1/2})_{\bullet i}||^2+0
\end{array}\\
\begin{array}{ll}
\phantom{\tr(\vH \vSigma\vH^\top)}&=\sum_{k\in \N_\ell}\left( ||(\vH\vSigma^{1/2})_{\bullet m_{k,1}}||^2+...+||(\vH\vSigma^{1/2})_{\bullet m_{k,1|\mathcal{M}_k|}}||^2\right)\\

&=\sum_{k\in \N_\ell}|| (\widetilde\vH\vSigma^{1/2})_{\bullet k}||^2\\
&=\sum_{k \in \N_\ell} (\widetilde \vH \vSigma \widetilde \vH^\top)_{kk}
\\&=\tr(\widetilde\vH \vSigma\widetilde\vH^\top).

\end{array}$\\

Here we used that if $\vH_{\bullet i}$ is a zero row, also $(\vH\vSigma^{1/2})_{\bullet i}\vT(\vX)$ is a zero row and that $\vH_{\bullet m_{k,j}}=\beta_{k,j} \cdot \vH_{\bullet \min(\mathcal M_k)}$. Together with the techniques from the first part, this leads to the result.\\ 
\end{proof}
\color{black}
From the proof of the corollary, it is clear that it is also possible to remove zero rows of $\vH$ without modifying the matrix regarding rows which are multiples of other rows. 
\\
Another operation which does not affect the value of the ATS is changing the order of the rows, while even elementary row operations (e.g., the addition of rows) in general affect the results.

\newpage
\bibliography{Literatur}

\end{document}